\documentclass[12pt,reqno]{amsart}
\usepackage{amssymb,amsthm,amsmath,amstext,amsxtra}
\usepackage{mathrsfs,bm}
\usepackage{fullpage}
\usepackage[all]{xy}
\usepackage{booktabs}
\usepackage{mathtools}
\usepackage{hyperref}
\hypersetup{colorlinks=true,urlcolor=blue,citecolor=blue,linkcolor=blue}
\usepackage{enumerate}
\usepackage{ stmaryrd }
\usepackage{enumitem}
\usepackage{comment}
\usepackage{multirow}
\usepackage{colonequals}
\usepackage{tikz}
\usepackage{marginnote}
\usepackage[export]{adjustbox}

\usepackage{tikz}
\usepackage{tikz-cd}
\usepackage{xcolor}

\numberwithin{equation}{section}

\newtheorem{example}{\text{Example}}[section]
\newtheorem{theorem}{Theorem}[section]

\newtheorem{conjecture}[theorem]{Conjecture}
\newtheorem{remark}[theorem]{Remark}
\newtheorem{definition}[theorem]{Definition}

\begin{document}

\title{Deligne's conjecture and mirror symmetry}

\author{Wenzhe Yang}
\address{SITP, Physics Department, Stanford University, CA, 94305}
\email{yangwz@stanford.edu}

\begin{abstract}
In this paper, we will study the connections between the mirror symmetry of Calabi-Yau threefolds and Deligne's conjecture on the special values of the $L$-functions of critical motives. Using the theory of mirror symmetry, we will develop a method to compute the Deligne's period for a Calabi-Yau threefold in the mirror family of a one-parameter mirror pair. We will give two examples to show how this method works, and we will express the Deligne's period in terms of the classical periods of the threeform. Using this method, we will compute the Deligne's period of a Calabi-Yau threefold studied in a recent paper by Candelas, de la Ossa, Elmi and van Straten. Based on their numerical results, we will explicitly show that this Calabi-Yau threefold satisfies Deligne's conjecture. A second purpose of this paper is to introduce the Deligne's conjecture to the physics community, and provide further evidence that there might exist interesting connections between physics and number theory.
\end{abstract}

\maketitle
\setcounter{tocdepth}{1}
\vspace{-13pt}
\tableofcontents
\vspace{-13pt}

\section{Introduction}
In the paper \cite{DeligneL}, Deligne formulates a profound conjecture about the relation between the special value of the $L$-function of an algebraic variety at a critical integral point and the classical periods of this variety. It provides an important generalization to the BSD conjecture for elliptic curves, while further prompts Beilinson to formulate a much more general conjecture about the special values of $L$-functions \cite{Nekovar, Schneider}. However, Deligne's conjecture is potentially extremely difficult to prove, and by far no proof is available in literature yet. Therefore it is very interesting to see whether researches from other areas, e.g. string theory and mirror symmetry, can provide any insights in the study of Deligne's conjecture, which is exactly the motivation of this paper. 

In fact, given a variety $X$, even the explicit computation of its Deligne's period is far from trivial. One important result of this paper is that for a Calabi-Yau threefold in a mirror family, mirror symmetry will provide all the geometric data needed in the computation of its Deligne's period. We will also give two important examples to show how this method works.

In order to streamline the discussion of Deligne's conjecture, we will resort to the concept of pure motives \cite{DeligneL}. But in this paper we will adopt the philosophy that pure motives can be understood through their classical realizations \cite{KimYang, Nekovar}. Given a smooth variety $X$ defined over $\mathbb{Q}$, a pure motive associated to $X$ will be denoted by $h^i(X)$, where $i$ is an integer. We will not attempt to explain what $h^i(X)$ is, or how to construct it. Instead, we will focus on its three important realizations:
\begin{enumerate}
\item The Betti realization $H^i_B(X)$, which is the singular cohomology group $H^i(X, \mathbb{Q})$ of the complex manifold structure of $X$ \cite{Nekovar}. It has a pure Hodge structure induced by the Hodge decomposition.

\item The de Rham realization $H^i_{\text{dR}}(X)$, which is constructed from the algebraic data of $X$, i.e. the algebraic forms on $X$ \cite{Voisin}.

\item The \'etale realization $H^i_{\text{\'et}}(X_{\overline{\mathbb{Q}}},\mathbb{Q}_\ell)$, which is the \'etale cohomology group of $X$ \cite{MilneEC}.

\end{enumerate}

The \'etale realization of $h^i(X)$ is a continuous representation of the absolute Galois group $\text{Gal}(\overline{\mathbb{Q}}/\mathbb{Q})$, which allows us to construct an $L$-function $L(h^i(X),s)$ \cite{Nekovar, Schneider}. The Tate twist of $h^i(X)$ by $\mathbb{Q}(n)$, $n \in \mathbb{Z}$ will be denoted by $h^i(X)(n)$, whose $L$-function satisfies 
\begin{equation}
L(h^i(X)(n),s)=L(h^i(X),n+s).
\end{equation}
Therefore, in order to study the special value of $L(h^i(X),s)$ at $s=n$, it is necessary and sufficient to study the special value of $L(h^i(X)(n),s)$ at $s=0$. The pure motive $h^i(X)(n)$ is critical if and only if its Hodge numbers satisfy the conditions given in the paper \cite{DeligneL}, which will be discussed in Section \ref{sec:Delignesconjecture}. In this paper, the most important example of a critical pure motive is $h^3(X)(2)$ when $X$ is a Calabi-Yau threefold defined over $\mathbb{Q}$. Deligne's conjecture claims that for a critical pure motive $h^i(X)(n)$, the special value $L(h^i(X)(n),0)$ is a rational multiple of its Deligne's period $c^+(h^i(X)(n))$.

The Deligne's period is constructed from the Betti realization and de Rham realization of $h^i(X)(n)$ \cite{DeligneL}. In order to explicitly compute the Deligne's period $c^+(h^i(X)(n))$, first we need to construct a rational basis for the Betti realization $H^i_B(X)(n)$. The complex conjugation defines an involution $F_\infty$ on $H^i_B(X)(n)$, and we also need to find the matrix of $F_\infty$ with respect to this rational basis. Next, we need to construct a rational basis for the de Rham realization $H^i_{\text{dR}}(X)(n)$, and explicitly write down the Hodge filtration of $H^i_{\text{dR}}(X)(n)$ with respect to this basis. In practice, it is very difficult to have all the previous data available, hence the computation of $c^+(h^i(X)(n))$ is certainly non-trivial. But this is exactly where mirror symmetry comes to the rescue! 

More explicitly, suppose we are given a one-parameter mirror pair of Calabi-Yau threefolds, where one-parameter means the Hodge number $h^{2,1}$ of the mirror threefold is $1$. We will focus on the case where the mirror threefold admits an algebraic deformation defined over $\mathbb{Q}$
\begin{equation} \label{eq:intromirrorfamily}
\pi:\mathscr{X} \rightarrow \mathbb{P}^1_{\mathbb{Q}}.
\end{equation}
In this paper, we will develop a method to compute the Deligne's period $c^+(h^3(\mathscr{X}_\varphi)(2))$ of a smooth rational fiber $\mathscr{X}_\varphi$ with $\varphi \in \mathbb{Q}$. Then we will give two important examples to illustrate how this method works. In order to verify that the critical motive $h^3(\mathscr{X}_\varphi)(2)$ satisfies Deligne's conjecture, we also need the special value of its $L$-function at $s=0$, i.e. $L(h^3(\mathscr{X}_\varphi)(2),0)$. However, generally it is very difficult to compute the special values of the $L$-functions of Calabi-Yau threefolds, even numerically.

In the paper \cite{Candelas}, the authors are able to find the $L$-function of a fiber $\mathscr{X}_{-1/7}$ over $\varphi=-1/7$ in the one-parameter mirror family called AESZ34 \cite{Duco, Candelas}. Using numerical method, they have found that the $L$-function $L(h^3(\mathscr{X}_{-1/7}),s)$ is given by
\begin{equation}
L(h^3(\mathscr{X}_{-1/7}),s)=L(f_2,s-1)L(f_4,s).
\end{equation}
Here $f_2$ is a weight-2 modular form for the modular group $\Gamma_0(14)$ that is labeled as \textbf{14.2.a.a} in LMFDB. While $f_4$ is a weight-4 modular form also for $\Gamma_0(14)$, which is labeled as \textbf{14.4.a.a} in LMFDB \cite{Candelas}. The authors have numerically computed the values of $L(f_2,1)$, $L(f_4,1)$ and $L(f_4,2)$. They have also numerically computed the values of the period matrix of the threeform at the point $\varphi=-1/7$, which can be expressed in terms of the values of $L(f_2,1)$, $L(f_4,1)$, $L(f_4,2)$ and $v^\perp$. Here $v^\perp$ is a number that is related to the period of the modular curve $X_0(14)$. The authors have speculated the connection between their results and Deligne's conjecture. Nevertheless they have not computed the Deligne's period of $h^3(\mathscr{X}_{-1/7})(2)$, hence they have not checked whether Deligne's conjecture is satisfied or not. 

Using the method developed in this paper, we are able to obtain an expression of the Deligne's period $c^+(h^3(\mathscr{X}_{-1/7})(2))$. Based on the numerical results in \cite{Candelas}, we have found that
\begin{equation}
c^+(h^3(\mathscr{X}_{-1/7})(2))=-\frac{2401}{32} L(f_2,1)L(f_4,2)=-\frac{2401}{32} L(h^3(\mathscr{X}_{-1/7})(2),0),
\end{equation} 
where the coefficient $-2401/32$ of course depends on the special rational bases chosen in our computation. Hence, we have verified that the critical pure motive $h^3(\mathscr{X}_{-1/7})(2)$ indeed satisfies Deligne's conjecture. We have also found that the period $c^-(h^3(\mathscr{X}_{-1/7}))$, also defined by Deligne in \cite{DeligneL}, satisfies the following equation
\begin{equation}
c^-(h^3(\mathscr{X}_{-1/7}))=\frac{1029}{32}\,\pi^{-3} \frac{L(f_4,1)L(f_2,1)}{v^\perp},
\end{equation}
while a detailed interpretation is presented in the paper \cite{Yang1}.

The outline of this paper is as follows. In Section \ref{sec:puremotive}, we will briefly review the theory of pure motives through the Betti, de Rham and \'etale realizations. In Section \ref{sec:Lfunction}, we will discuss the $L$-functions associated to pure motives. In Section \ref{sec:Delignesconjecture}, we will discuss the construction of Deligne's periods and introduce Deligne's profound conjecture. In Section \ref{sec:mirrorsymmetry}, we will briefly review the theory of the mirror symmetry of Calabi-Yau threefolds. In Section \ref{sec:computationdeligneperiods}, we will develop a method to compute the Deligne's periods of Calabi-Yau threefolds in a one-parameter mirror family. In Section \ref{sec:explicitDeligneperiods}, we will use the method in Section \ref{sec:computationdeligneperiods} to explicitly compute the Deligne's periods of two important examples. In Section \ref{sec:exampleAESZ34}, we will compute the Deligne's period of the special Calabi-Yau threefold $\mathscr{X}_{-1/7}$ in the mirror family AESZ34, and verify that it satisfies Deligne's conjecture. In Section \ref{sec:conclusion} we will summarize the results of this paper and discuss several open questions.

\section{The pure motives} \label{sec:puremotive}

In this section, we will briefly introduce the theory of pure motives through their classical realizations \cite{Nekovar}. The language of pure motives, even though very abstract, can greatly simplify the studies of arithmetic geometry and number theory \cite{KNY, Nekovar, Schneider}. In this section, we will assume $X$ is a smooth projective variety defined over $\mathbb{Q}$. Let $M$ be the following pure motive associated to $X$
\begin{equation} \label{eq:puremotiveM}
M:=h^i(X)(n),~i, n \in \mathbb{Z}.
\end{equation}
We now explain the meaning of $M$ through its three classical realizations:
\begin{enumerate}
\item The Betti realization. The $\mathbb{C}$-valued points (classical points) of $X$, denoted by $X(\mathbb{C})$, form a smooth projective complex manifold. More concretely, $X(\mathbb{C})$ is just the complex manifold associated to $X$ in the usual sense. The Betti realization of $M$, denoted by $M_{\text{B}}$, is just the following singular cohomology group
\begin{equation}
M_{\text{B}}:=H^i\big(X(\mathbb{C}),~\mathbb{Q}(n) \big)=H^i\big(X(\mathbb{C}),~\mathbb{Q} \big) \otimes \mathbb{Q}(n),
\end{equation}
where $\mathbb{Q}(n)$ means the rational vector space $(2\pi i)^n \, \mathbb{Q}$ \cite{Nekovar}. Recall that $\mathbb{Q}(n)$ admits a pure Hodge structure of Hodge type $(-n,-n)$ \cite{PetersSteenbrink}. The Hodge decomposition
\begin{equation} \label{eq:BettiHodgeD}
M_{\text{B}} \otimes_{\mathbb{Q}} \mathbb{C}= \oplus_{p+q=w}\, H^{p,q}
\end{equation}
defines a pure Hodge structure on $M_{\text{B}}$ with weight $w:=i-2\,n$. Together with this pure Hodge structure, the Betti realization induces the Hodge realization of $M_{\text{B}}$. By definition, the Hodge number $h^{p,q}$ is
\begin{equation}
h^{p,q}:=\text{dim}_{\mathbb{C}}\,H^{p,q}.
\end{equation}
The complex conjugation $c \in \text{Gal}(\mathbb{C}/\mathbb{R})$ defines an action on the points of $X(\mathbb{C})$, which further induces an involution $c^*$ on $M_{B}$. Let $F_\infty$ be the involution on $M_{B}$ induced by the action of $c$ on both the points of $X(\mathbb{C})$ and the coefficient ring $\mathbb{Q}(n)$. Then the conjugate-linear involution $F_\infty \otimes c$ preserves the Hodge decomposition \ref{eq:BettiHodgeD}, i.e. it sends $H^{p,q}$ to $H^{p,q}$.

\item The de Rham realization. On the variety $X$, there exists a complex of sheaves of algebraic differential forms \cite{Hartshorne}
\begin{equation} 
\Omega_X^*: 0 \rightarrow \mathcal{O}_X \xrightarrow{d} \Omega_{X}^1 \xrightarrow{d} \cdots \xrightarrow{d}  \Omega_{X}^{\text{dim}(X)} \rightarrow 0.
\end{equation}
In order to define a `reasonable' cohomology theory, first we need to choose an injective resolution $\Omega_{X}^* \rightarrow I^* $ in the abelian category that consists of the complexes of sheaves on $X$. Then the algebraic de Rham cohomology of $X$ is defined by \cite{Voisin}
\begin{equation}
\mathbb{H}^i(X_{\text{Zar}},\Omega^*_{X}):=H^i(\Gamma(X,I^*)),
\end{equation}
which is also called the hypercohomology of $\Omega_{X}^*$. Here $X_{\text{Zar}}$ means the Zariski topology on $X$. The de Rham realization of $M$, denoted by $M_{\text{dR}}$, is just the hypercohomology of the shifted complex of sheaves $\Omega^*_{X}[n]$
\begin{equation}
M_{\text{dR}}:=\mathbb{H}^i(X_{\text{Zar}},\Omega^*_{X}[n]),~\text{where}~(\Omega^*_{X}[n])^l=\Omega^{l+n}_{X},
\end{equation}
which is in fact a finite dimensional vector space over $\mathbb{Q}$ \cite{Voisin}. Furthermore, $M_{\text{dR}}$ has a decreasing filtration $F^p M_{\text{dR}}$ given by
\begin{equation}
F^p M_{\text{dR}}:=\mathbb{H}^i(X_{\text{Zar}}, F^p\Omega_{X}^*[n]),
\end{equation}
where the complex $F^p\Omega_{X}^*[n]$ is of the form
\begin{equation}
F^p\, \Omega_{X}^*[n]: 0 \rightarrow \cdots \rightarrow 0 \rightarrow \Omega_{X}^{p+n} \xrightarrow{d}  \Omega_{X}^{p+1+n} \xrightarrow{d} \cdots \xrightarrow{d}  \Omega_{X}^{\text{dim}\,X} \rightarrow 0.
\end{equation}

\item The $\ell$-adic realization. Suppose $\ell$ is a prime number, then the $\ell$-adic cohomology of $X$ is by definition the inverse limit
\begin{equation}
H^i_{\text{\'et}}(X_{\overline{\mathbb{Q}}},\mathbb{Q}_\ell):=\varprojlim_n H^i((X \times_\mathbb{Q} \overline{\mathbb{Q}})_{\text{\'et}}, \mathbb{Z}/\ell^n \mathbb{Z}) \otimes_{\mathbb{Z}_\ell} \mathbb{Q}_{\ell},
\end{equation}
where $(X \times_\mathbb{Q} \overline{\mathbb{Q}})_{\text{\'et}}$ is the \'etale topology on the $\overline{\mathbb{Q}}$-variety $X_{\overline{\mathbb{Q}}}:=X \times_\mathbb{Q} \overline{\mathbb{Q}}$ and $\mathbb{Z}/\ell^n \mathbb{Z}$ means the constant \'etale torsion sheaf on $(X \times_\mathbb{Q} \overline{\mathbb{Q}})_{\text{\'et}}$. The $\ell$-adic cyclotomic character $\mathbb{Q}_\ell(1)$ is the inverse limit
\begin{equation}
\mathbb{Q}_\ell(1):= \varprojlim_n \mu_{\ell^n}(\overline{\mathbb{Q}}) \otimes_{\mathbb{Z}_{\ell}} \mathbb{Q}_\ell,
\end{equation}
where $\mu_{\ell^n}(\overline{\mathbb{Q}})$ consists of the $\ell^n$-th root of unity \cite{Taylor}. Let $\mathbb{Q}_\ell(n)$ be the $n$-fold tensor product $\mathbb{Q}_\ell(1)^{\otimes n}$, which is a continuous representation of $\text{Gal}(\overline{\mathbb{Q}}/\mathbb{Q})$ \cite{Taylor}. The $\ell$-adic realization of $M$, denoted by $M_\ell$, is given by
\begin{equation}
M_\ell:=H^i_{\text{\'et}}(X_{\overline{\mathbb{Q}}},\mathbb{Q}_\ell) \otimes_{\mathbb{Q}_\ell} \mathbb{Q}_{\ell}(n),
\end{equation}  
which is also a continuous representation of $\text{Gal}(\overline{\mathbb{Q}}/\mathbb{Q})$  \cite{MilneEC}. 
\end{enumerate}

There exist standard comparison isomorphisms between the three realizations \cite{Nekovar}:
\begin{enumerate}
\item There is an isomorphism $I_\infty$ between the Betti realization and the de Rham realization
\begin{equation} \label{eq:bettiderhamComparison}
I_\infty:M_{B} \otimes_{\mathbb{Q}} \mathbb{C} \rightarrow M_{\text{dR}} \otimes \mathbb{C},
\end{equation}
which sends $\oplus_{k \geq p} H^{k,w-k}$ to $F^pM_{\text{dR}} \otimes \mathbb{C}$. This comparison isomorphism $I_\infty$ sends the involution $F_ \infty \otimes c$ on $M_{B} \otimes_{\mathbb{Q}} \mathbb{C}$ to the involution $1 \otimes c$ on $M_{\text{dR}} \otimes \mathbb{C}$. This property will be crucial when we compute the Deligne's period for a Calabi-Yau threefold in a one-parameter mirror family. 

\item Suppose $\overline{\infty}: \overline{\mathbb{Q}} \hookrightarrow \mathbb{C}$ is an embedding of $\overline{\mathbb{Q}} $ into $\mathbb{C}$, then there exists an isomorphism $I_{\ell,\overline{\infty}}$ between the Betti realization and the $\ell$-adic realization 
\begin{equation}
I_{\ell,\overline{\infty}}:M_{B} \otimes_{\mathbb{Q}} \mathbb{Q}_\ell \rightarrow M_{\ell},
\end{equation}
which depends on the choice of $\overline{\infty}$ up to an isomorphism. Moreover, $I_{\ell,\overline{\infty}}$ sends the involution $F_\infty \otimes 1$ on $M_{B} \otimes_{\mathbb{Q}} \mathbb{Q}_\ell$ to the involution $c$ on $M_{\ell}$.
\end{enumerate}

\noindent 
The two comparison isomorphisms immediately imply that
\begin{equation}
\text{dim}_{\mathbb{Q}}(M_{B})=\text{dim}_\mathbb{Q}(M_{\text{dR}})=\text{dim}_{\mathbb{Q}_\ell}(M_\ell),
\end{equation}
and the common dimension is denoted by $\text{dim}(M)$, which is called the rank of $M$.
\begin{example}
The Tate motive $\mathbb{Q}(1)$ is by definition the dual of the Lefschetz motive $h^2(\mathbb{P}^1_{\mathbb{Q}})$, whose classical realizations are:
\begin{enumerate}
\item $\mathbb{Q}(1)_{\text{B}}=(2 \pi i)\, \mathbb{Q}$, which admits a pure Hodge structure of type $(-1,-1)$.

\item $\mathbb{Q}(1)_{\text{dR}}=\mathbb{Q}$, with Hodge filtration given by $F^0=0$ and $F^{-1}=\mathbb{Q}$.

\item $\mathbb{Q}(1)_{\ell}=\mathbb{Q}_\ell(1)$.
\end{enumerate}
The Tate motive $\mathbb{Q}(n)$ is the $n$-fold tensor product $\mathbb{Q}(1)^{\otimes n}$ \cite{Nekovar, Schneider}.
\end{example}

\noindent The Tate twist of the pure motive $M$ by $\mathbb{Q}(m)$ is by definition the tensor product
\begin{equation}
M(m):=M \otimes \mathbb{Q}(m).
\end{equation}
So the pure motive $M$ in the formula \ref{eq:puremotiveM} can also be expressed as
\begin{equation}
M=h^i(X) \otimes \mathbb{Q}(n).
\end{equation}
There exist a Poincar\'e duality and a hard Lefschetz theorem for each of the three classical realizations of pure motives, which are compatible with each other under the previous standard comparison isomorphisms. Therefore the dual of $M$ is given by \cite{Nekovar}
\begin{equation} \label{eq:dualpuremotive}
M^\vee=h^i(X)^\vee(-n)=h^{2\text{dim}X-i}(X)(\text{dim}X-n)=h^i(X)(i-n)=M(w);~w=i-2n
\end{equation} 

Intuitively, we can say the pure motives of the form \ref{eq:puremotiveM} encode all of the cohomological information of the smooth projective variety $X$. While the theory of pure motives can be viewed as a universal cohomology theory that bridges algebraic geometry and number theory \cite{KimYang, Nekovar}.

\section{The \texorpdfstring{$L$}{L}-function of a pure motive} \label{sec:Lfunction}

In this section, we will introduce the $L$-function associated to a pure motive. The reader who is familiar with the construction can skip this section completely.

Given a pure motive $M=h^i(X)(n)$, its $\ell$-adic realization $M_\ell$ is a continuous representation of $\text{Gal}(\overline{\mathbb{Q}}/\mathbb{Q})$ \cite{Nekovar,Taylor}. The inertia group $I_p$ for a prime number $p$ is a subgroup of $\text{Gal}(\overline{\mathbb{Q}}/\mathbb{Q})$, whose definition can be found in the book \cite{Nekovar}. The representation $M_\ell$ is said to be unramified at $p$ if the action of $I_p$ on $M_\ell$ is trivial. If so, the geometric Frobenius $\text{Fr}_p$ has a well-defined action on $M_\ell$ \cite{SerreLF, Taylor}. Since $X$ is a smooth projective variety defined over $\mathbb{Q}$,  $M_{\ell}$ is pure of weight $w=i-2n$. Here `pure' means that there exists a finite set $\mathcal{S}$ consisting of prime numbers such that if $p$ satisfies $p \notin \mathcal{S}$ and $p \nmid \ell$, $M_\ell$ is unramified at $p$.

On the other hand, for every prime number $p \neq \ell$, let $M_\ell^{I_p}$ be the subspace of $M_\ell$ that is fixed under the action of $I_p$. Then the geometric Frobenius has a well-defined action on $M_\ell^{I_p}$, and its characteristic polynomial is given by
\begin{equation}
P_p(M,T)=\text{det}\big( 1-T\,\text{Fr}_p \vert M_\ell^{I_p} \big),~ \ell \neq p.
\end{equation}
From Deligne's proof of Weil conjectures \cite{DeligneWeil}, if $X$ has good reduction at $p$, then we have:
\begin{enumerate}
\item $P_p(M,T)$ is an integral polynomial and it is independent of the choice of $\ell$.

\item $P_p(M,T)$ has a factorization of the form
\begin{equation}
P_p(M,T)= \prod_{j=1}^{\text{dim}(M)}(1-\alpha_j \, T),
\end{equation}
where $\alpha_j$ is an algebraic integer with $|\alpha_j|=p^{w/2}$ for every $j$.
\end{enumerate}
The variety $X$ only has bad reduction at finitely many primes. Serre has a conjecture about the properties of $P_p(M,T)$ at a bad prime of $X$ \cite{SerreL}.
\begin{conjecture} \label{eq:Serresconjecture}
Given an arbitrary prime number $p$, $P_p(M,T)$ is an integral polynomial which does not depend on the choice of $\ell$. It can be factorized into
\begin{equation}
P_p(M,T)= \prod_{j=1}^{\text{dim}( M_\ell^{I_p})}(1-\alpha_j \, T),
\end{equation}
where for every $j$, $\alpha_j$ is an algebraic integer with absolute value
\begin{equation}
|\alpha_j|=p^{w_j/2},~0 \leq w_j \leq w.
\end{equation}

\end{conjecture}

By definition, the local $L$-factor of $M$ at $p$ is
\begin{equation} \label{eq:defnlocalLfactor}
L_p(M,s):=\frac{1}{P_p(M,p^{-s})},
\end{equation}
and the $L$-function of $M$ is defined by
\begin{equation} \label{eq:defnLfunction}
L(M,s):= \prod_p L_p(M,s),
\end{equation}
where the infinite product is over all the prime numbers. The local $L$-factor $L_p(M,s)$ satisfies the following properties \cite{Nekovar}
\begin{equation} \label{eq:transportLfunction}
L_p(M(m),s)=L_p(M,m+s),~~L_p(M_1 \oplus M_2,s)=L_p(M_1,s)L_p(M_2,s);
\end{equation} 
hence the $L$-function of $M$ satisfies similar properties, i.e.
\begin{equation} \label{eq:translationLfunction}
L(M(m),s)=L(M,m+s),~L(M_1 \oplus M_2,s)=L(M_1,s)L(M_2,s).
\end{equation}
Deligne's proof of Weil conjectures and the Conjecture \ref{eq:Serresconjecture} by Serre imply that $L(M,s)$ converges absolutely when $\text{Re}(s) > w/2+1$, thus $L(M,s)$ is a nowhere vanishing holomorphic function in this region. However, the existence of a meromorphic extension of $L(M,s)$ to the complex plane is still a conjecture \cite{DeligneL,Schneider}.

\begin{conjecture}
There exists a meromorphic extension of $L(M,s)$ to the entire complex plane $\mathbb{C}$. When the weight $w$ is odd, this extension is globally holomorphic, while when $w$ is even, the only possible pole of this extension is at $s=w/2+1$. Furthermore, if $s=w/2+1$ is not a pole, then the special value $L(M,w/2+1)$ is non-zero.
\end{conjecture}

The archimedean prime of $\mathbb{Z}$ is the natural embedding of $\mathbb{Q}$ into $\mathbb{C}$, which will be denoted by $\infty$ \cite{Nekovar}. There is also a local $L$-factor $L_\infty(M,s)$ associated to the archimedean prime, but its construction is more involved \cite{DeligneL, Neukirch}. For simplicity, let us define the Gamma factors $\Gamma_{\mathbb{R}}(s)$ and $\Gamma_{\mathbb{C}}(s)$ by
\begin{equation}
\begin{aligned}
\Gamma_{\mathbb{R}}(s):&=\pi^{-s/2}\cdot \Gamma(s/2),\\
\Gamma_{\mathbb{C}}(s):&=\Gamma_{\mathbb{R}}(s)\cdot \Gamma_{\mathbb{R}}(s+1)=2\cdot (2 \pi)^{-s}\cdot \Gamma(s),
\end{aligned}
\end{equation}
where $\Gamma(s)$ is the Gamma function. The local $L$-factor $L_\infty(M,s)$ only depends on the pure Hodge structure on $M_{\text{B}} \otimes_{\mathbb{Q}} \mathbb{R}$, and its construction is carefully discussed in \cite{SerreL} and Section 5.2 of \cite{DeligneL}:
\begin{enumerate}
\item If the weight $w$ of $M$ is odd, then $L_\infty(M,s)$ is defined by
\begin{equation}
L_\infty(M,s)= \prod_{p < q} \Gamma_{\mathbb{C}}(s-p)^{h^{p,q}} ~.
\end{equation}

\item When the weight $w$ of $M$ is even, the subspace $H^{w/2,w/2}$ in the Hodge decomposition can be further decomposed according to the eigenvalues of the involution $F_\infty$ into
\begin{equation}
H^{w/2,w/2}=H^{w/2,+} \oplus H^{w/2,-}; 
\end{equation}
where $H^{w/2,+}$ and $H^{w/2,-}$ are determined by the condition
\begin{equation}
F_\infty|_{H^{w/2,+}}=(-1)^{w/2},~F_\infty|_{H^{w/2,-}}=(-1)^{w/2+1}.
\end{equation}
Then $L_\infty(M,s)$ is defined to be
\begin{equation}
\begin{aligned}
L_\infty(M,s)= \prod_{p < q} \Gamma_{\mathbb{C}}(s-p)^{h^{p,q}} &\cdot \Gamma_{\mathbb{R}}(s-w/2)^{\text{dim}H^{w/2,+}} \\
&\cdot \Gamma_{\mathbb{R}}(s-w/2+1)^{\text{dim}H^{w/2,-}}.
\end{aligned}
\end{equation}
\end{enumerate} 
The local $L$-factor $L_\infty(M,s)$ also satisfies the following properties \cite{DeligneL}
\begin{equation} \label{eq:translationInfinity}
L_\infty(M(m),s)=L_\infty(M,m+s),~L_\infty (M_1 \oplus M_2,s)=L_\infty(M_1,s)\cdot L_\infty(M_2,s).
\end{equation} 
The full $L$-function of $M$ is given by
\begin{equation}
\Lambda(M,s)=L(M,s)\cdot L_{\infty}(M,s).
\end{equation}

\begin{conjecture}
$\Lambda(M,s)$ satisfies the functional equation \cite{DeligneL, SerreL}
\begin{equation} \label{eq:functionalequation}
\Lambda(M,s)=\varepsilon(M,s)\,\Lambda(M^{\vee},1-s),
\end{equation} 
where $\varepsilon(M,s)$ is of the form $a \cdot b^s$ with $a$ and $b$ being non-zero complex numbers. From the formula \ref{eq:dualpuremotive}, this functional equation is equivalent to
\begin{equation} 
\Lambda(M,s)=\varepsilon(M,s)\,\Lambda(M,w+1-s).
\end{equation} 
\end{conjecture}

The study of the special values of the $L$-function $L(M,s)$ at integral points $s \in \mathbb{Z}$ has been a central theme in modern number theory. These special values of $L$-functions are closely related to the geometric information of the pure motive $M$.

\section{Deligne's conjecture} \label{sec:Delignesconjecture}

In this section, we will introduce Deligne's conjecture on the special values of the $L$-functions of critical pure motives, and we will follow the original paper \cite{DeligneL} closely. First, let us explain what is a critical pure motive.
\begin{definition}
For a pure motive $M$, an integer $n$ is called critical if neither $L_\infty(M,s)$ nor $L_\infty(M^\vee,1-s)$ has a pole at $s=n$.
\end{definition}
Deligne's conjecture concerns the special value of $L(M,s)$ at a critical integer $s=n$. Since we always have the freedom to twist $M$ by $\mathbb{Q}(n)$, so the formulas \ref{eq:translationLfunction} and \ref{eq:translationInfinity} imply that there is no loss of generality if we assume $n=0$ from the beginning. A pure motive $M$ is said to be critical if $n=0$ is critical for $M$. In fact, whether a pure motive is critical or not is completely determined by its Hodge numbers \cite{DeligneL}. More precisely, $M$ is critical if and only if the following conditions are satisfied:
\begin{enumerate}
\item For every pair $(p,q)$ of integers such that $p\neq q$ and $h^{p,q}\neq 0$, we have $p \leq -1, ~q\geq 0$ or $p \geq 0,~ q \leq -1$.

\item If the weight $w$ of $M$ is even, then the action of $F_\infty$ on $H^{w/2,w/2}$ is 1 if $w < 0$ and $-1$ if $w \geq 0$.
\end{enumerate}
For the purpose of this paper, we will only consider the case where the weight $w$ of $M$ is odd.

Let $M_B^+$ and $M_B^-$ be the subspaces of the Betti realization$M_B$ of $M$ defined by the conditions
\begin{equation}
F_\infty |_{M_B^+} = 1, ~F_\infty |_{M_B^-} = -1;
\end{equation}
and let $d^{+}(M)$ and $d^{-}(M)$ be their dimensions respectively
\begin{equation}
d^{+}(M)=\text{dim}_\mathbb{Q} \, M_B^+,~d^{-}(M)=\text{dim}_\mathbb{Q}\,M_B^-.
\end{equation}
Since the involution $F_\infty$ maps $H^{p,q}$ onto $H^{q,p}$ and vice versa, thus if the weight $w$ of $M$ is odd, we must have
\begin{equation}
d^{+}(M)=d^{-}(M)=\frac{1}{2}\text{dim}_\mathbb{Q}\,(M_B).
\end{equation}
On the other hand, let $F^+$ and $F^-$ be the linear subspaces occurring in the Hodge filtration $F^pM_{\text{dR}}$ such that
\begin{equation}
\text{dim}_\mathbb{Q}\,F^+ =\text{dim}_\mathbb{Q}\,M_B^+,~\text{dim}_\mathbb{Q}\,F^- =\text{dim}_\mathbb{Q}\,M_B^-.
\end{equation}
More explicitly, via the comparison isomorphism between the Betti and de Rham realizations, $F^+ \otimes \mathbb{C}$ corresponds to 
\begin{equation}
\oplus_{p > q} H^{p,q}(M_B) ~\text{with}~p+q=w.
\end{equation}
If the weight $w$ of $M$ is odd, then we will have $F^-=F^+$. Next, let us define $M_{\text{dR}}^\pm$ by
\begin{equation}
M_{\text{dR}}^+:=M_{\text{dR}}/F^-,~M_{\text{dR}}^-:=M_{\text{dR}}/F^+,
\end{equation}
while if $w$ is odd, we must have $M_{\text{dR}}^+=M_{\text{dR}}^-$. 

The comparison isomorphism $I_\infty$ \ref{eq:bettiderhamComparison} and the natural projection map $M_{\text{dR}} \rightarrow M_{\text{dR}}^+$ induce a composition of maps denoted by $I^+_\infty$
\begin{equation}
I^+_\infty:M_B^+ \otimes \mathbb{C} \hookrightarrow M_B \otimes \mathbb{C} \xrightarrow{I_\infty} M_{\text{dR}} \otimes \mathbb{C} \rightarrow M_{\text{dR}}^+ \otimes \mathbb{C}.
\end{equation}
Since $F_\infty$ swaps $H^{p,q}$ and $H^{q,p}$, the homomorphism $I^+_\infty$ is in fact an isomorphism \cite{DeligneL}. Now choose a rational basis for $M_B^+$ and a rational basis for $M_{\text{dR}}^+$. Then the determinant of $I^+_\infty$ with respect to the two rational bases can be computed, which by definition is the Deligne's period of $M$
\begin{equation}
c^+(M)=\text{det}(I^+_\infty).
\end{equation}
Notice that the Deligne's period $c^+(M)$ is only well-defined up to a nonzero rational multiple! By exactly the same construction, there is another isomorphism $I^-_\infty$ of the form
\begin{equation}
I^-_\infty:M_B^- \otimes \mathbb{C} \rightarrow M_{\text{dR}}^- \otimes \mathbb{C}.
\end{equation}
Similarly, a rational basis for $M_B^-$ and a rational basis for $M_{\text{dR}}^-$ allow us to define another period $c^-(M)$ by
\begin{equation}
c^-(M)=\text{det}(I^-_\infty),
\end{equation}
which is also well-defined up to a nonzero rational multiple. The definition of the Deligne's period $c^+(M)$ (or $c^-(M)$) means that it can be expressed in terms of the classical periods. Notice that the dual of $M_{\text{dR}}^+$ is the subspace $F^+M_{\text{dR}}^\vee$ of the dual $M_{\text{dR}}^\vee$. If we choose a basis $\{\omega_i \}$ for $F^+M_{\text{dR}}^\vee$ and a basis $\{\rho_i \}$ for $M_B^+$, then the matrix of $I^+_\infty$ with respect to these two bases is given by $\langle\omega_i,\rho_j \rangle$. Here the pairing $\langle ,\rangle$ is induced by the Poincar\'e duality. The Deligne's period $c^+(M)$ is just the determinant of this matrix
\begin{equation}
c^+(M)=\text{det}(\langle\omega_i,\rho_j \rangle).
\end{equation}
Similarly, we can also express $c^-(M)$ in terms of classical periods. However, generally it is not easy to explicitly compute the Deligne's period, as it can be very difficult to find the rational bases for $M_B^+$ and $F^+M_{\text{dR}}^\vee$. But for a Calabi-Yau threefold in a one-parameter mirror family, we will see that the mirror symmetry will provide all the data needed in the computation of its Deligne's period.

From the paper \cite{DeligneL}, if $M$ is critical, then $s=0$ is not a pole of $L(M,s)$, i.e. the special value of $L(M,s)$ at $s=0$ is a finite number. Now, we are ready to state the conjecture of Deligne, which is about the relation between $L(M, 0)$ and $c^+(M)$ when $M$ is critical.

\vspace*{0.1in}

\textbf{Deligne's conjecture}: \textit{If the pure motive $M$ is critical, then $L(M,0)$ is a rational multiple of $c^+(M)$.}

\vspace*{0.1in}

A proof of Deligne's conjecture is still not available in literature, therefore it is very important to provide interesting examples for it, which might shed lights on the nature of the conjecture itself.

\section{The Mirror symmetry of Calabi-Yau threefolds} \label{sec:mirrorsymmetry}

In this section, we will briefly review the mirror symmetry of Calabi-Yau threefolds, and we will focus on the one-parameter mirror pairs \cite{PhilipXenia, CoxKatz, MarkGross, KimYang}. Given a mirror pair $(X^\vee,X)$ of Calabi-Yau threefolds, one-parameter means that the Hodge numbers of $X^\vee$ and $X$ satisfy
\begin{equation}
h^{1,1}(X^\vee)=h^{2,1}(X)=1.
\end{equation}

\subsection{The Picard-Fuchs equation}

For the purpose of this paper, we will assume that the mirror threefold $X$ has an algebraic deformation defined over $\mathbb{Q}$
\begin{equation} \label{eq:mirrorfamily}
\pi:\mathscr{X} \rightarrow \mathbb{P}^1_{\mathbb{Q}},
\end{equation}
where the coordinate of the base $\mathbb{P}^1_{\mathbb{Q}}$ is denoted by $\varphi$. From now on, $X$ will also mean the underlying differential manifold structure of a smooth fiber of the family \ref{eq:mirrorfamily}. We will further assume that for each smooth fiber $\mathscr{X}_\varphi$, there exists a nowhere-vanishing algebraic threeform $\Omega_\varphi$ that varies algebraically with respect to $\varphi$. Moreover, as a form on $\mathscr{X}$, $\Omega$ is defined over $\mathbb{Q}$. Hence for a rational point $\varphi$, $\Omega_\varphi$ is also defined over $\mathbb{Q}$ \cite{CoxKatz, MarkGross, KimYang}. From the Griffiths transversality, the threeform $\Omega_\varphi$ satisfies a fourth-order Picard-Fuchs equation
\begin{equation} \label{eq:PFequation}
\mathscr{L}\,\Omega_\varphi=0,
\end{equation}
where $\mathscr{L}$ is a differential operator with polynomial coefficients $R_i(\varphi) \in \mathbb{Q}[\varphi]$
\begin{equation} \label{eq:PicardFuchsOperator}
\mathscr{L}=R_4(\varphi) \,\vartheta^4+R_3(\varphi)\, \vartheta^3 +R_2(\varphi)\, \vartheta^2+ R_1(\varphi) \, \vartheta^1+R_0(\varphi), ~ \vartheta=\varphi \,\frac{d}{d\varphi}.
\end{equation}

The Picard-Fuchs operator $\mathscr{L}$ has finitely many regular singularities, and a singularity is called the large complex structure limit if the monodromy at it is maximally unipotent \cite{CoxKatz, MarkGross}. In this paper, we will assume that the Picard-Fuchs operator $\mathcal{L}$ has a large complex structure limit at $\varphi=0$. More concretely, there exists a small neighborhood $\Delta$ of $\varphi=0$, on which the Picard-Fuchs equation \ref{eq:PFequation} has four canonical solutions of the form 
\begin{equation} \label{eq:PeriodsCan}
\begin{aligned}
\varpi_0 &= f_0,  \\
\varpi_1 &=\frac{1}{2\pi i}\left(f_0 \log \varphi+f_1\right), \\
\varpi_2 &=\frac{1}{(2\pi i)^2}\left( f_0 \log^2 \varphi +2\, f_1 \log \varphi + f_2\right), \\
\varpi_3 &=\frac{1}{(2 \pi i)^3} \left( f_0 \log^3 \varphi +3 \, f_1 \log^2 \varphi +3\, f_2 \log \varphi +f_3 \right),
\end{aligned} 
\end{equation}
where $\{f_j\}_{j=0}^3$ are power series that converge on $\Delta$. If we further impose the  conditions
\begin{equation} \label{eq:boundarycondition}
f_0(0)=1,~f_1(0)=f_2(0)=f_3(0)=0,
\end{equation}
then the four canonical solutions \ref{eq:PeriodsCan} become unique. The canonical period vector $\varpi$ is the following column vector
\begin{equation}
\varpi:=(\varpi_0,\,\varpi_1,\,\varpi_2,\,\varpi_3)^\top.
\end{equation}
\begin{remark}
In this paper, the multi-valued holomorphic function $\log \varphi$ is chosen to satisfy
\begin{equation}
\log (1)=0,~\log (-1)=\pi i.
\end{equation}
\end{remark}

The Poincar\'e duality implies that there exists a unimodular skew symmetric pairing on the homology group $H_3(X,\mathbb{Z})$ (modulo torsion), which allows us to choose an integral symplectic basis $\{A_0,A_1,B_0,B_1\}$ that satisfies the following intersection pairing \cite{PhilipXenia,CoxKatz,MarkGross}
\begin{equation}
A_a \cdot A_b=0,~~B_a \cdot B_b=0,~~A_a \cdot B_b= \delta_{ab}.
\end{equation}
Suppose the dual of this basis is $\{\alpha^0,\alpha^1,\beta^0,\beta^1\}$, i.e. we have
\begin{equation}
\alpha^a (A_b)=\delta_{ab}, ~ \beta^a(B_b)=\,\delta_{ab},~\alpha^a (B_b)=\beta^a(A_b)=0,
\end{equation}
then it forms a basis for $H^3(X,\mathbb{Z})$ (modulo torsion). From the Poincar\'e duality, we have 
\begin{equation} \label{eq:cupproductHatcher}
\int_X \alpha^a \smile \beta^b=\delta_{ab},~\int_X \alpha^a \smile \alpha^b=0,~\int_X \beta^a \smile \beta^b=0,
\end{equation}
where $\alpha^a \smile \beta^b$ means the cup product between $\alpha^a$ and $\beta^b$, etc \cite{Hatcher}.

\begin{remark}
The torsion of the homology or cohomology groups are irrelevant to this paper, hence they will be ignored.
\end{remark}

The integration of the threeform $\Omega_\varphi$ on the symplectic basis $\{A_a,B_a\}_{a=0}^1$ gives us the integral periods
\begin{equation} \label{eq:IntegralPeriodDefinition}
z_a(\varphi)= \int_{A_a} \Omega_\varphi,~\mathcal{G}_b(\varphi)=\int_{B_b} \Omega_\varphi,
\end{equation}
which are multi-valued holomorphic functions \cite{PhilipXenia,CoxKatz,MarkGross}. Now let us define the integral period vector $\amalg(\varphi)$ to be the column vector
\begin{equation} \label{eq:integralperiodvector}
\amalg(\varphi):=(\mathcal{G}_0(\varphi),\mathcal{G}_1(\varphi),z_0(\varphi),z_1(\varphi))^\top.
\end{equation}
Since the integral period vector $\amalg$ forms another basis for the solution space of the Picard-Fuchs equation \ref{eq:PFequation}, there exists a transformation matrix $S \in \text{GL}(4,\mathbb{C})$ such that
\begin{equation} \label{eq:PiSomega}
\amalg=S\cdot \varpi.
\end{equation}
The transformation matrix $S$ is crucial in this paper, and it can be evaluated by mirror symmetry \cite{PhilipXenia, KimYang}. For later convenience, let us also define the row vector $\beta$ by
\begin{equation}
\beta:=(\beta^0,\beta^1,\alpha^0,\alpha^1).
\end{equation}
Under the comparison isomorphism, the threeform $\Omega_\varphi$ has an expansion of the form
\begin{equation} \label{eq:expansionOmega}
\Omega_\varphi=\beta \cdot \amalg(\varphi)=\mathcal{G}_0(\varphi) \, \beta^0 +\mathcal{G}_1(\varphi) \,\beta^1+z_0(\varphi)\, \alpha^0+z_1(\varphi)\, \alpha^1.
\end{equation}

\subsection{The prepotential}

For a one-parameter mirror pair $(X^\vee,X)$, the complexified K\"ahler moduli space $\mathscr{M}_K(X^\vee)$ of $X^\vee$ has a very simple description \cite{MarkGross,KimYang}
\begin{equation}
\mathscr{M}_K(X^\vee)=(\mathbb{R}+i\, \mathbb{R}_{> 0})/\mathbb{Z}=\mathbb{H}/\mathbb{Z},
\end{equation}
where $\mathbb{H}$ is the upper half plane of $\mathbb{C}$. Now let $e$ be a basis of $H^2(X^\vee,\mathbb{Z})$ (modulo torsion) which also lies in the K\"ahler cone of $X^\vee$, then every point of $\mathscr{M}_K(X^\vee)$ can be represented by $e\,t, \,t \in \mathbb{H}$ \cite{MarkGross}. While $e\,t$ is equivalent to $e\,(t+1)$ under the quotient by $\mathbb{Z}$. In physics literature, $t$ is called the flat coordinate of $\mathscr{M}_K(X^\vee)$ \cite{PhilipXenia, CoxKatz, MarkGross}. The prepotential $\mathcal{F}$ admits an expansion near $t=i\, \infty$ given by \cite{PhilipXenia,CoxKatz}
\begin{equation} \label{eq:Prepotential}
\mathcal{F}=-\frac{1}{6}\, Y_{111}\, t^3 -\frac{1}{2}\, Y_{011}\,t^2-\frac{1}{2}\,Y_{001}\, t-\frac{1}{6}\,Y_{000}+\mathcal{F}^{\text{np}},
\end{equation}
where $\mathcal{F}^{\text{np}}$ is the non-perturbative instanton correction. Moreover, $\mathcal{F}^{\text{np}}$ is invariant under the translation $t \rightarrow t+1$ and it is also exponentially small when $t \rightarrow i\, \infty$, i.e. it admits a series expansion in $\exp 2\pi i\,t$
\begin{equation}
\mathcal{F}^{\text{np}}=\sum_{n=1}^{\infty} a_n \exp 2 \pi i \,nt.
\end{equation}
The coefficient $Y_{111}$ in \ref{eq:Prepotential} is the topological intersection number given by \cite{PhilipXenia,CoxKatz, MarkGross}
\begin{equation} \label{eq:y111}
Y_{111}=\int _{X^\vee} e\wedge e \wedge e,
\end{equation}
which is a positive integer. The computations of the coefficients $Y_{011}$ and $Y_{001}$ are more tricky, but in this paper, we will only need the fact that they are rational numbers \cite{KimYang}. In all examples of mirror pairs, $Y_{000}$ is always of the form \cite{PhilipXenia}
\begin{equation} \label{eq:PhysicistsY000}
Y_{000}=-3\, \chi (X^\vee)\, \frac{\zeta(3)}{(2 \pi i)^3},
\end{equation}
where $\chi(X^\vee)$ is the Euler characteristic of $X^\vee$. A detailed study of the occurrence of $\zeta(3)$ from the motivic point of view is presented in the paper \cite{KimYang}.

\subsection{The mirror symmetry}

In all examples of one-parameter mirror pairs, there exists an integral symplectic basis $\{A_0,A_1,B_0,B_1\}$ of $H^3(X,\mathbb{Z})$ such that \cite{CoxKatz,MarkGross}
\begin{equation} \label{eq:zivarpii01}
z_j(\varphi)=\lambda (2 \pi i)^3 \,\varpi_j(\varphi),~j=0,1;~ \lambda \in \mathbb{Q}^\times.
\end{equation}
Let us denote the quotient $\varpi_1/\varpi_0$ by $t_c$, which is of the form
\begin{equation}
t_c=\frac{z_1}{z_0}=\frac{\varpi_1}{\varpi_0}=\frac{1}{2 \pi i}\,\log \varphi+\frac{f_1(\varphi)}{f_0(\varphi)}.
\end{equation}
Under the action of the monodromy, i.e. $\log \varphi \rightarrow \log \varphi+2 \pi i$, $t_c$ transforms in the way
\begin{equation}
t_c \rightarrow t_c +1.
\end{equation} 

The mirror map for a one-parameter mirror pair is induced by the identification of the coordinate $t$ on the K\"ahler side and the coordinate $t_c$ on the complex side
\begin{equation} \label{eq:Mirrormap}
t \equiv t_c.
\end{equation}
Hence from now on, we will use the notations $t$ and $t_c$ interchangeably. The normalization of the integral period vector $\amalg$ \ref{eq:integralperiodvector} is denoted by $\amalg_{\text{A}}$
\begin{equation}
\amalg_{\text{A}}=(\mathcal{G}_0/z_0,\mathcal{G}_1/z_0,\,1,\,z_1/z_0)^\top.
\end{equation} 
On the K\"ahler side, the period vector $\Pi$ is determined by the prepotential $\mathcal{F}$ \ref{eq:Prepotential} \cite{PhilipXenia,KimYang}
\begin{align}
\Pi =(\mathcal{F}_0,\mathcal{F}_1, 1,t)^\top,~\text{with}~\mathcal{F}_0= 2 \,\mathcal{F}-t \,\frac{\partial  \mathcal{F}}{\partial t},~ \mathcal{F}_1=\frac{\partial \mathcal{F}}{\partial t}.
\end{align}
The mirror symmetry claims that under the mirror map \ref{eq:Mirrormap}, we have the following equation
\begin{equation}
\Pi=\amalg_{\text{A}}.
\end{equation}
Now we are ready to compute the transformation matrix $S$ in the formula \ref{eq:PiSomega}. Near the large complex structure limit $\varphi=0$, the conditions in the formula \ref{eq:boundarycondition} imply 
\begin{equation}
t=\frac{1}{2 \pi i} \,\log \varphi+ \mathcal{O}(\varphi).
\end{equation}
Therefore, under the mirror map \ref{eq:Mirrormap}, $\varphi=0$ on the complex side corresponds to $t= i\, \infty$ on the K\"ahler side \cite{PhilipXenia, CoxKatz, KimYang}. In the limit $t \rightarrow i\, \infty$, the leading parts of $\amalg_{\text{A}}$ and  $\varpi$ are given by
\begin{align} \label{eq:Pilimit}
\amalg_{\text{A}} \equiv \Pi \sim
\begin{pmatrix}
\frac{1}{6}\, Y_{111}\, t^3 -\frac{1}{2}\, Y_{001}\, t-\frac{1}{3}\, Y_{000} \\
-\frac{1}{2}\,Y_{111} \,t^2 - Y_{011}\,t-\frac{1}{2} \,Y_{001} \\
1 \\
t \\
\end{pmatrix}
,~ \varpi \sim
\begin{pmatrix}
1 \\
t \\
t^2 \\
t^3 \\
\end{pmatrix},
\end{align}
from which the transformation matrix $S$ can be easily evaluated \cite{KimYang}
\begin{equation} \label{eq:smatrixrepn}
S\,=\lambda(2 \pi i)^3 \,
\left(
\begin{array}{cccc}
 -\frac{1}{3}\, Y_{000} & -\frac{1}{2} \,Y_{001} & 0 & \frac{1}{6}\, Y_{111} \\
 -\frac{1}{2} \,Y_{001} & -\,Y_{011} & -\frac{1}{2} \,Y_{111} & 0 \\
 1 & 0 & 0 & 0 \\
 0 & 1 & 0 & 0 \\
\end{array}
\right),~\lambda \in \mathbb{Q}^\times.
\end{equation}

\section{The computations of the Deligne's period in mirror symmetry} \label{sec:computationdeligneperiods}

In this section, we will apply the results in the previous section to the computation of the Deligne's period for a smooth fiber $\mathscr{X}_\varphi, \varphi \in \mathbb{Q}$ in a one-parameter mirror family.

Given such a smooth fiber $\mathscr{X}_\varphi, \varphi \in \mathbb{Q}$, the algebraic de Rham cohomology group $H^3_{\text{dR}}(\mathscr{X}_\varphi)$ is a four dimensional vector space that has an explicit basis. Since the threeform $\Omega_\varphi$ is nowhere vanishing on $\mathscr{X}_\varphi$, hence it lies in the subspace $ F^3 H^3_{\text{dR}}(\mathscr{X}_\varphi)$. But $F^3 H^3_{\text{dR}}(\mathscr{X}_\varphi) $ is one dimensional, therefore $\Omega_\varphi$ must form a basis for it. Let the derivative of $\Omega_\varphi$ with respect to $\varphi$ be denoted by $\Omega'_\varphi$. From the Griffiths transversality, $\Omega'_\varphi$ lies in $F^2 H^3_{\text{dR}}(\mathscr{X}_\varphi) $. But the dimension of $F^2 H^3_{\text{dR}}(\mathscr{X}_\varphi) $ is two, therefore $\{\Omega_\varphi,\Omega'_\varphi \}$ must form a basis for it \cite{KimYang}. Similarly, we deduce that $\{\Omega_\varphi,\Omega'_\varphi,\Omega''_\varphi \}$ forms a basis for $F^1 H^3_{\text{dR}}(\mathscr{X}_\varphi) $ and $\{\Omega_\varphi,\Omega'_\varphi,\Omega''_\varphi, \Omega'''_\varphi \}$ forms a basis for $F^0 H^3_{\text{dR}}(\mathscr{X}_\varphi) $ \cite{KimYang}. Put everything together, we have 
\begin{equation} \label{eq:HodgefiltrationMS}
\begin{aligned}
F^3(H^3_{\text{dR}}(\mathscr{X}_\varphi))&=\langle \Omega_\varphi \rangle, \\
F^2(H^3_{\text{dR}}(\mathscr{X}_\varphi))&=\langle \Omega_\varphi,\Omega'_\varphi \rangle, \\
F^1(H^3_{\text{dR}}(\mathscr{X}_\varphi))&=\langle \Omega_\varphi,\Omega'_\varphi,\Omega''_\varphi \rangle, \\
F^0(H^3_{\text{dR}}(\mathscr{X}_\varphi))&=\langle \Omega_\varphi,\Omega'_\varphi,\Omega''_\varphi ,\Omega'''_\varphi \rangle;
\end{aligned}
\end{equation}
where $\langle \Omega_\varphi \rangle$ means the vector space spanned by $\Omega_\varphi$ over $\mathbb{Q}$, etc. Under the comparison isomorphism, the derivatives of $\Omega_\varphi$ have expansions given by
\begin{equation} \label{eq:derivativesOmega}
\Omega^{(n)}_\varphi=\beta \cdot S \cdot \varpi^{(n)},~n=0,1,2,3;
\end{equation}
where we have used the formulas \ref{eq:expansionOmega} and \ref{eq:PiSomega} \cite{KimYang}. By definition, the Wronskian $W$ of the period vector $\varpi$ is
\begin{equation}
W\,=
\left(
\begin{array}{cccc}
 \varpi_0 & \varpi'_0 & \varpi''_0 & \varpi'''_0 \\
 \varpi_1 & \varpi'_1 & \varpi''_1 & \varpi'''_1 \\
 \varpi_2 & \varpi'_2 & \varpi''_2 & \varpi'''_2 \\
 \varpi_3 & \varpi'_3 & \varpi''_3 & \varpi'''_3 \\
\end{array}
\right),
\end{equation}
the determinant of which does not vanish at a smooth point $\varphi$ \cite{Yang}. Under the comparison isomorphism, the rational basis $(\Omega_\varphi,\Omega'_\varphi,\Omega''_\varphi, \Omega'''_\varphi )$ of $ H^3_{\text{dR}}(\mathscr{X}_\varphi) \otimes \mathbb{C} $ is mapped to the basis $\beta \cdot S \cdot W$ of $H^3(X,\mathbb{Q}) \otimes \mathbb{C}$. Recall that $X$ also means the underlying differential manifold structure of a smooth fiber $\mathscr{X}_\varphi$.

Given a smooth point $\varphi \in \mathbb{Q}$, the action of the involution $F_\infty$ on the Betti cohomology $H^3(X,\mathbb{Q})$ can be computed explicitly. The key property needed in the computation is that under the comparison isomorphism, the involution $F_\infty \otimes c$ on $H^3(X,\mathbb{Q}) \otimes \mathbb{C}$ corresponds to the involution $1 \otimes c$ on $ H^3_{\text{dR}}(\mathscr{X}_\varphi) \otimes \mathbb{C} $. From this property, we immediately deduce that
\begin{equation}
\beta \cdot S \cdot W=\beta \cdot F_\infty \cdot \overline{S} \cdot \overline{W},
\end{equation}
where $\overline{S}$ (resp. $\overline{W}$) means the complex conjugation of the matrix $S$ (resp. $W$). Thus the matrix of $F_\infty$ with respect to the basis $\beta$ of $H^3(X,\mathbb{Q})$ is given by
\begin{equation} \label{eq:complexconjugationinvolution}
F_\infty=S \cdot W \cdot \overline{W}^{-1} \cdot \overline{S}^{-1}.
\end{equation}
In fact, the involution $F_\infty$ is also defined on the integral cohomology group $H^3(X,\mathbb{Z})$
\begin{equation}
F_\infty: H^3(X,\mathbb{Z}) \rightarrow H^3(X,\mathbb{Z}),
\end{equation} 
therefore with respect to the integral symplectic basis $\beta$ of $H^3(X,\mathbb{Z})$, $F_\infty$ is an  integral matrix. Let $\mathcal{R}$ be the set of real-valued regular singularities of the Picard-Fuchs operators $\mathscr{L}$ \ref{eq:PicardFuchsOperator}, i.e.
\begin{equation}
\mathcal{R}=\{ \varphi \in \mathbb{R}:\mathscr{L}~\text{is singular at }~\varphi \}.
\end{equation}
The entries of the matrix
\begin{equation} \label{eq:finfinitysmoothrepn}
S \cdot W \cdot \overline{W}^{-1} \cdot \overline{S}^{-1}
\end{equation}
are smooth functions on $\varphi \in \mathbb{R}- \mathcal{R}$, whose values at a rational point are integral. Since the rational points are dense in  $\mathbb{R}- \mathcal{R}$, we immediately deduce that the matrix \ref{eq:finfinitysmoothrepn} is locally constant on $\mathbb{R}- \mathcal{R}$, hence it is constant in every open interval of $\mathbb{R}- \mathcal{R}$.

From the Poincar\'e duality, the dual of $ H^3_{\text{dR}}(\mathscr{X}_\varphi)$ is given by \cite{Hartshorne, Voisin}
\begin{equation}
H^3_{\text{dR}}(\mathscr{X}_\varphi)^\vee=H^3_{\text{dR}}(\mathscr{X}_\varphi) \otimes \mathbb{Q}(3).
\end{equation} 
From Section \ref{sec:Delignesconjecture}, its subspace $F^+(H^3_{\text{dR}}(\mathscr{X}_\varphi) \otimes \mathbb{Q}(3))$ is
\begin{equation}
F^+(H^3_{\text{dR}}(\mathscr{X}_\varphi) \otimes \mathbb{Q}(3))=F^{-1}(H^3_{\text{dR}}(\mathscr{X}_\varphi) \otimes \mathbb{Q}(3)),
\end{equation}
which is the two dimensional vector space spanned by $\Omega_\varphi$ and $\Omega'_\varphi$. Suppose the subspace of $H^3(X,\mathbb{Q})$ on which $F_\infty$ acts as 1 has a basis $(\gamma^+_0,\gamma^+_1)$, then from Section \ref{sec:Delignesconjecture}, the period $c^+(h^3(\mathscr{X}_\varphi))$ is given by
\begin{equation}
c^+(h^3(\mathscr{X}_\varphi))=  \det 
\left(
\begin{array}{cc}
 \frac{1}{(2 \pi i)^3} \int_{X} \Omega_\varphi \smile \gamma^+_0  &  \frac{1}{(2 \pi i)^3} \int_{X} \Omega_\varphi \smile \gamma^+_1  \\
 \frac{1}{(2 \pi i)^3} \int_{X} \Omega'_\varphi \smile \gamma^+_0  &  \frac{1}{(2 \pi i)^3} \int_{X} \Omega'_\varphi \smile \gamma^+_1  \\
\end{array}
\right).
\end{equation}
Here the cup products can be computed by the formulas \ref{eq:cupproductHatcher} and \ref{eq:derivativesOmega}. Notice that the additional factor $(2 \pi i)^{-3}$ comes from the fact that the dual of $H^3(X,\mathbb{Q})$ is $H^3(X,\mathbb{Q}) \otimes \mathbb{Q}(3)$, and the pairing
\begin{equation}
H^3(X,\mathbb{Q}) \times ( H^3(X,\mathbb{Q}) \otimes \mathbb{Q}(3) ) \rightarrow \mathbb{Q}
\end{equation}
is given by
\begin{equation}
\langle \phi_1, \phi_2 \rangle=\frac{1}{(2 \pi i)^3} \int_X \phi_1 \smile \phi_2; \phi_1 \in H^3(X,\mathbb{Q}), \phi_2 \in H^3(X,\mathbb{Q}) \otimes \mathbb{Q}(3).
\end{equation}
The period $c^-(h^3(\mathscr{X}_\varphi))$ can be computed similarly. Suppose the subspace of $H^3(X,\mathbb{Q})$ on which $F_\infty$ acts as $-1$ has a basis $(\gamma^-_0,\gamma^-_1)$, then $c^-(h^3(\mathscr{X}_\varphi))$ is given by
\begin{equation}
c^-(h^3(\mathscr{X}_\varphi))= \det 
\left(
\begin{array}{cc}
 \frac{1}{(2 \pi i)^3} \int_{X} \Omega_\varphi \smile \gamma^-_0  &  \frac{1}{(2 \pi i)^3} \int_{X} \Omega_\varphi \smile \gamma^-_1  \\
 \frac{1}{(2 \pi i)^3} \int_{X} \Omega'_\varphi \smile \gamma^-_0  &  \frac{1}{(2 \pi i)^3} \int_{X} \Omega'_\varphi \smile \gamma^-_1  \\
\end{array}
\right).
\end{equation}
The upshot is that the Deligne's periods $c^\pm(h^3(\mathscr{X}_\varphi))$ can be explicitly expressed in terms of the values of the canonical periods $\varpi_i$ and its derivatives.

Since the Hodge numbers of a smooth fiber $\mathscr{X}_\varphi$ in the mirror family \ref{eq:mirrorfamily} are
\begin{equation}
h^{3,0}=h^{2,1}=h^{1,2}=h^{0,3}=1,
\end{equation}
thus from Section \ref{sec:Delignesconjecture}, the pure motive $h^3(\mathscr{X}_\varphi)(n)$ is critical if and only if when $n=2$ \cite{DeligneL}. A basis for the Betti cohomology group $H^3(X,\mathbb{Q}) \otimes \mathbb{Q}(2)$ is given by $(2 \pi i)^2 \beta$. Since $(2 \pi i)^2$ is a real number, the subspace of $H^3(X,\mathbb{Q})\otimes \mathbb{Q}(2)$ on which $F_\infty$ acts as 1 has a basis $((2 \pi i)^2\gamma^+_0,(2 \pi i)^2\gamma^+_1)$. Therefore, we immediately have \cite{DeligneL}
\begin{equation}
c^+(h^3(\mathscr{X}_\varphi)(2))=(2 \pi i)^4 c^+(h^3(\mathscr{X}_\varphi)).
\end{equation}
Similarly, we also have
\begin{equation}
c^-(h^3(\mathscr{X}_\varphi)(2))=(2 \pi i)^4 c^-(h^3(\mathscr{X}_\varphi)).
\end{equation}
\begin{remark}
Since the values of $\varpi_i$ and its derivatives at a point $\varphi$ can be computed numerically to a very high precision, hence $c^\pm(h^3(\mathscr{X}_\varphi)(2))$ can also be evaluated numerically.
\end{remark}

\section{Two important examples } \label{sec:explicitDeligneperiods}

In this section, we will apply the method developed in Section \ref{sec:computationdeligneperiods} to compute the Deligne's periods for two important examples. Suppose $\varphi_{-1}$ and $\varphi_1$ are two real singularities of the Picard-Fuchs operator $\mathscr{L}$ \ref{eq:PicardFuchsOperator} such that
\begin{equation}
\varphi_{-1} <0 < \varphi_1.
\end{equation}
Suppose further that $\mathscr{L}$ does not have any other singularities in the interval $(\varphi_{-1},0)$ or $(0,\varphi_1)$, i.e. $\varphi_{-1}$ is the largest negative singularity of $\mathscr{L}$ and $\varphi_1$ is the smallest positive singularity of $\mathscr{L}$. Recall that the singularity $\varphi=0$ is the large complex structure limit of  $\mathscr{L}$. In this section, we will explicitly compute the Deligne's period for a rational fiber $\mathscr{X}_\varphi,\varphi \in \mathbb{Q}$ such that
\begin{equation}
0< \varphi < \varphi_1 ~\text{or}~\varphi_{-1}< \varphi <0.
\end{equation} 
Let us first look at the case where $0< \varphi < \varphi_1$.

\subsection{The first case}

Let $\varphi \in \mathbb{Q}$ be a small positive number such that the power series $f_i$ in \ref{eq:PeriodsCan} converges at it. Since $f_i$ is a power series with rational coefficients, i.e. it lies in $\mathbb{Q}[[\varphi]]$, we deduce that $\varpi^{(n)}_0(\varphi)$ and $\varpi_2^{(n)}(\varphi)$ and   are real numbers, while $\varpi_1^{(n)}(\varphi)$ and $\varpi_3^{(n)}(\varphi)$ are purely imaginary numbers. Hence we have $W=V\cdot \overline{W}$, where $V$ is the matrix
\begin{equation}
V\,=
\left(
\begin{array}{cccc}
 1 & 0 & 0 & 0 \\
 0 & -1& 0 & 0 \\
 0 & 0 & 1 & 0 \\
 0 & 0 & 0 & -1 \\
\end{array}
\right).
\end{equation}
Then from the formula \ref{eq:complexconjugationinvolution}, we deduce that the matrix of $F_\infty$ with respect to the basis $\beta$ of $H^3(X,\mathbb{Q})$ is given by
\begin{equation} \label{eq:involutionpositivephi}
F_\infty=S \cdot V \cdot \overline{S}^{-1}=
\left(
\begin{array}{cccc}
 1 & 0 & 0 & 0 \\
 0 & -1 & 0 & -2 Y_{011} \\
 0 & 0 & -1 & 0 \\
 0 & 0 & 0 & 1 \\
\end{array}
\right).
\end{equation}
But from Section \ref{sec:computationdeligneperiods}, we know $F_\infty$ is constant in the open interval $(0,\varphi_1)$. Hence, for every rational point $\varphi \in (0,\varphi_1)$, the matrix of $F_\infty$ is given by the formula \ref{eq:involutionpositivephi}. The two linearly independent eigenvectors of $F_\infty$ associated to the eigenvalue 1 are
\begin{equation} \label{eq:eigenv1positive}
(0,-Y_{011},0,1)^\top~\text{and}~(1,0,0,0)^\top.
\end{equation}
While the two linearly independent eigenvectors associated to the eigenvalue $-1$ are 
\begin{equation} \label{eq:eigenvm1positive}
(0,0,1,0)^\top~\text{and}~(0, 1, 0, 0)^\top.
\end{equation}
Therefore, the subspace of $H^3(X,\mathbb{Q})$ on which $F_\infty$ acts as 1 has a basis
\begin{equation}
\beta^0,~\alpha^1-Y_{011} \beta^1,
\end{equation}
and the subspace of $H^3(X,\mathbb{Q})$ on which $F_\infty$ acts as $-1$ has a basis
\begin{equation}
\alpha^0,~\beta^1.
\end{equation}

The cup product between $\Omega_\varphi$ and the two eigenvectors in the formula \ref{eq:eigenv1positive} are given by 
\begin{equation}
\begin{aligned}
\int_X\Omega_\varphi \smile \beta^0 &=\lambda (2 \pi i)^3 \varpi_0,\\
\int_X \Omega_\varphi \smile (\alpha^1-Y_{011} \beta^1) &=\lambda (2 \pi i)^3\left(\frac{1}{2} Y_{001} \varpi_0+\frac{1}{2}Y_{111} \varpi_2 \right).
\end{aligned}
\end{equation}
Hence the Deligne's period $c^+(h^3(\mathscr{X}_\varphi))$ is
\begin{equation}
c^+(h^3(\mathscr{X}_\varphi))=\frac{1}{2}\lambda^2Y_{111}  (\varpi_0\varpi'_2-\varpi_2 \varpi'_0).
\end{equation}
where $\lambda$ is the non-zero rational constant in the formula \ref{eq:smatrixrepn} and $Y_{111}$ is the topological intersection number in the formula \ref{eq:y111}. From the definition of Deligne's period, we have the freedom to rescale $c^+(h^3(\mathscr{X}_\varphi))$ by a non-zero rational number. So, we can simply let $c^+(h^3(\mathscr{X}_\varphi))$ be
\begin{equation}
c^+(h^3(\mathscr{X}_\varphi))=\varpi_0\varpi'_2-\varpi_2 \varpi'_0.
\end{equation}
Similarly, $c^-(h^3(\mathscr{X}_\varphi))$ is given by
\begin{equation}
\begin{aligned}
c^-(h^3(\mathscr{X}_\varphi))=\frac{1}{6}\lambda^2Y_{111}\Big( (\frac{2Y_{000}}{Y_{111}}\varpi_0-\varpi_3) \varpi'_1-(\frac{2Y_{000}}{Y_{111}}\varpi'_0-\varpi'_3) \varpi_1 \Big).
\end{aligned}
\end{equation}
We can again throw away the overall nonzero rational constant in the expression, and simply let $c^-(h^3(\mathscr{X}_\varphi))$ be
\begin{equation}
c^-(h^3(\mathscr{X}_\varphi))=(\frac{2Y_{000}}{Y_{111}}\varpi_0-\varpi_3) \varpi'_1-(\frac{2Y_{000}}{Y_{111}}\varpi'_0-\varpi'_3) \varpi_1.
\end{equation}

\begin{remark}
The numbers $Y_{111}$ and $Y_{000}$ are determined by the topological data, which are independent of the choice of a symplectic basis $\beta$. While the numbers $Y_{011}$ and $Y_{001}$ do depend on the choice of $\beta$, but they do not occur in the expressions of $c^\pm(h^3(\mathscr{X}_\varphi))$.
\end{remark}

It is very interesting to notice that both $c^+(h^3(\mathscr{X}_\varphi))$ and $c^-(h^3(\mathscr{X}_\varphi))$ are expressed as the values of multi-valued holomorphic functions at $\varphi$. More precisely, the exist two multi-valued holomorphic functions whose values at a rational point $\varphi \in (0,\varphi_1)$ give us $c^+(h^3(\mathscr{X}_\varphi))$ and $c^-(h^3(\mathscr{X}_\varphi))$.

\subsection{The second case} \label{sec:negativecase}

Now suppose $\varphi \in \mathbb{Q}$ is a small negative number such that the power series $f_i$ in \ref{eq:PeriodsCan} converges at it. Since $f_i$ lies in $\mathbb{Q}[[\varphi]]$, the values of $f_i$ and its derivatives at $\varphi$ are real numbers. Recall that by our choice, $\log(-1)=\pi i$ and $\log 1=0$, so we have
\begin{equation}
\log \varphi=\log(-\varphi)+\pi i; ~\varphi \in \mathbb{R}~\text{and}~\varphi<0,
\end{equation}
where $\log(-\varphi)$ is a real number. Thus under complex conjugation, we deduce that
\begin{equation}
\overline{ \log \varphi }=\log \varphi-2 \pi i.
\end{equation}
Therefore, we obtain $W=V\cdot \overline{W}$, where $V$ is given by
\begin{equation}
V\,=
\left(
\begin{array}{cccc}
 1 & 0 & 0 & 0 \\
 1 & -1 & 0 & 0 \\
 1 & -2 & 1 & 0 \\
 1 & -3 & 3 & -1 \\
\end{array}
\right).
\end{equation}
Then the formula \ref{eq:complexconjugationinvolution} tells us that the matrix of $F_\infty$ with respect to the basis $\beta$ of $H^3(X,\mathbb{Q})$ is given by
\begin{equation} \label{eq:involutionnegativephi}
F_\infty=S \cdot V \cdot \overline{S}^{-1}=
\left(
\begin{array}{cccc}
 1 & 1 & Y_{001}-\frac{1}{6} Y_{111} & Y_{011}+\frac{1}{2} Y_{111} \\
 0 & -1 & Y_{011}+\frac{1}{2} Y_{111} & -2 Y_{011}-Y_{111} \\
 0 & 0 & -1 & 0 \\
 0 & 0 & -1 & 1 \\
\end{array}
\right).
\end{equation}
But from Section \ref{sec:computationdeligneperiods}, we know $F_\infty$ is constant in the open interval $ (\varphi_{-1},0)$. Hence, for every rational point $\varphi \in (\varphi_{-1},0)$, the matrix of $F_\infty$ is given by the formula \ref{eq:involutionnegativephi}. The two eigenvectors of $F_\infty$ associated to the eigenvalue 1 are
\begin{equation}\label{eq:eigenv1negative}
(0,-Y_{011}-\frac{1}{2} Y_{111},0,1)^\top~\text{and}~(1,0,0,0)^\top.
\end{equation}
While the two eigenvectors associated to the eigenvalue $-1$ are 
\begin{equation}\label{eq:eigenvm1negative}
(-Y_{001}-\frac{1}{2} Y_{011}-\frac{1}{12} Y_{111},0,2,1)^\top~\text{and}~(-\frac{1}{2},1,0,0)^\top.
\end{equation}
Therefore, the subspace of $H^3(X,\mathbb{Q})$ on which $F_\infty$ acts as 1 has a basis
\begin{equation}
\alpha^1+\left(-Y_{011}-\frac{1}{2} Y_{111} \right) \beta^1,~\beta^0,
\end{equation}
and the subspace of $H^3(X,\mathbb{Q})$ on which $F_\infty$ acts as $-1$ has a basis
\begin{equation}
2 \alpha^0+\alpha^1+ \left( -Y_{001}-\frac{1}{2} Y_{011}-\frac{1}{12} Y_{111}\right) \beta^0,~-\frac{1}{2} \beta^0+\beta^1.
\end{equation}
Similarly, $c^\pm(h^3(\mathscr{X}_\varphi))$ are given by
\begin{equation}
\begin{aligned}
&c^+(h^3(\mathscr{X}_\varphi)) = \varpi_0(\varpi'_2-\varpi'_1) -\varpi'_0(\varpi_2-\varpi_1), \\
&c^-(h^3(\mathscr{X}_\varphi)) =2\left(\varpi_3-\frac{3}{2} \varpi_2 \right)\varpi'_0-2\left(\varpi'_3-\frac{3}{2} \varpi'_2 \right)\varpi_0 + \\
& \left(\frac{8Y_{000}-Y_{111}}{Y_{111}}\varpi_0+6\varpi_2-4\varpi_3 \right) \varpi'_1-\left(\frac{8Y_{000}-Y_{111}}{Y_{111}}\varpi'_0+6\varpi'_2-4\varpi'_3 \right) \varpi_1.\\
\end{aligned}
\end{equation}
Again, we have found that there exist two multi-valued holomorphic functions whose values at a rational point $\varphi \in (\varphi_{-1},0)$ give us $c^+(h^3(\mathscr{X}_\varphi))$ and $c^-(h^3(\mathscr{X}_\varphi))$.

\begin{remark}
Notice that for a general smooth point $\varphi \in \mathbb{Q}$, $F_\infty$ can be evaluated numerically, which allows us to numerically compute $c^\pm(h^3(\mathscr{X}_\varphi))$.
\end{remark}

\section{An example for Deligne's Conjecture} \label{sec:exampleAESZ34}

In this section, we will use the method developed in the previous sections to compute the Deligne's period of a Calabi-Yau threefold that has been studied in the paper \cite{Candelas}. Based on their numerical results, we will explicitly verify that this Calabi-Yau threefold satisfies Deligne's conjecture.

More explicitly, in the paper \cite{Candelas}, the authors have studied the one-parameter mirror pair $(X^\vee,X)$ of Calabi-Yau threefolds that is called AESZ34 \cite{Duco, Meyer}. The mirror threefold $X$ has an algebraic deformation of the form
\begin{equation} \label{eq:aeszfamilies}
\pi:\mathscr{X} \rightarrow \mathbb{P}^1_{\mathbb{Q}}.
\end{equation}
The zeta functions of the smooth fiber $\mathscr{X}_{-1/7}$ over $\varphi=-1/7$ have been numerically computed for small prime numbers, from which the authors are able to find the $L$-function of the pure motive $h^3(\mathscr{X}_{-1/7})$. The numerical values of the canonical periods $\varpi_i$ (and their derivatives) at $\varphi=-1/7$ have also been computed by them to a very high precision, and they are able to express these values in terms of the special values of $L$-functions. They have speculated the connections between  their numerical results and Deligne's conjecture. But they have not computed the Deligne's period $c^+(h^3(\mathscr{X}_{-1/7})(2))$ for the critical motive $h^3(\mathscr{X}_{-1/7})(2)$, therefore Deligne's conjecture has not been numerically verified. In this section, we will use the method developed in Section \ref{sec:computationdeligneperiods} and Section \ref{sec:explicitDeligneperiods} to compute $c^+(h^3(\mathscr{X}_{-1/7})(2))$, then we will explicitly verify that Deligne's conjecture is satisfied by the critical motive $h^3(\mathscr{X}_{-1/7})(2)$. We will also numerically compute the period $c^-(h^3(\mathscr{X}_{-1/7}))$ and look at its properties.

\subsection{An overview of the mirror pair AESZ34}

First, let us review the results of the paper \cite{Candelas} that will be needed in this paper, while the readers are referred to it for more details. The Hodge diamond of the mirror threefold $X$ of the mirror pair AESZ34 is of the form \cite{Candelas}
\begin{center} 
\begin{tabular}{ c c c c c c c }
 &  &  & 1 &  &  &  \\ 
 &  & 0&   & 0&  &  \\   
 & 0&  & 9 &  & 0&   \\  
1&  & 1 &  & 1 & & 1 \\ 
 & 0&  & 9 &  & 0&   \\ 
 &  & 0&   & 0&  &  \\   
 &  &  & 1 &  &  &  \\
\end{tabular}.
\end{center}
The Picard-Fuchs equation of the algebraic deformation \ref{eq:aeszfamilies} of $X$, i.e. the mirror family, is
\begin{equation} \label{eq:aesz34PF}
\begin{aligned}
\mathcal{D}=& \theta^4 -\varphi(35 \theta^4+70 \theta^3+63 \theta^2+28 \theta +5) + \varphi^2 (\theta+1)^2(259 \theta^2+518 \theta +285)\\
                       &-225 \varphi^3 (\theta+1)^2(\theta+2)^2,~ \theta=\varphi \frac{d}{d \varphi}.
\end{aligned}
\end{equation}
This Picard-Fuchs operator $\mathcal{D}$ has five regular singularities
\begin{equation}
\varphi=0, 1/25, 1/9, 1, \infty,
\end{equation}
while $\varphi=0$ is the large complex structure limit. The canonical period $\varpi_0$ of $\mathcal{D}$ is  \cite{Duco}
\begin{equation}
\varpi_0=1+\sum_{n=1}^\infty a_n \varphi^n;~~a_n=\sum_{i+j+k+l+m=n}\left(\frac{n! }{i!j!k!l!m!} \right)^2.
\end{equation}
The numbers that occur in the perturbative part of the prepotential $\mathcal{F}$ \ref{eq:Prepotential} have been computed in \cite{Candelas}
\begin{equation} \label{eq:prepotentialcoeff}
Y_{111}=24,~Y_{011}=0,~Y_{001}=-2,~Y_{000}=48 \frac{\zeta(3)}{(2 \pi i)^3}.
\end{equation}

The zeta functions of the pure motive $h^3(\mathscr{X}_{-1/7})$ for small prime numbers have been numerically computed. At a good prime number $p$, the zeta function of $h^3(\mathscr{X}_{-1/7})$ has a factorization of the form
\begin{equation}
(1-a_p (pT)+p(pT)^2)(1-b_p T+p^3T^2).
\end{equation}
Here $a_p$ is the $p$-th coefficient of the $q$-expansion of a weight-2 modular form $f_2$ for the modular group $\Gamma_0(14)$, which is labeled as \textbf{14.2.a.a} in LMFDB. While $b_p$ is the $p$-th coefficient of the $q$-expansion of a weight-4 modular form $f_4$ also for $\Gamma_0(14)$, which is labeled as \textbf{14.4.a.a} in LMFDB. Notice that this property has only been numerically checked by them for small prime numbers \cite{Candelas}. Hence the $L$-function of the pure motive $h^3(\mathscr{X}_{-1/7})$ should be
\begin{equation}
L(h^3(\mathscr{X}_{-1/7}),s)=L(f_2,s-1)L(f_4,s).
\end{equation}
In particular, the special value $L(h^3(\mathscr{X}_{-1/7}),2)$ is just $L(f_2,1)L(f_4,2)$. In the paper \cite{Candelas}, both $L(f_2,1)$ and $L(f_4,2)$ have been numerically computed to a very high precision
\begin{equation}
\begin{aligned}
L(f_2,1)&=0.33022365934448053902826194612283487754045234078189 \cdots ,\\
L(f_4,2)&=0.91930674266912115653914356907939249680895763199044  \cdots.
\end{aligned}
\end{equation}

The power series expansions of the canonical periods $\varpi_i$ do not converge at $\varphi=-1/7$, nevertheless their values can be computed to a very high precision by numerically solving the Picard-Fuchs equation \ref{eq:aesz34PF}. In \cite{Candelas}, the numerical values of the canonical periods $\varpi_i$ (and their derivatives) at $\varphi=-1/7$ have been computed, which can be expressed in terms of the special values $L(f_2,1)$, $L(f_4,1)$, $L(f_4,2)$ and $v^\perp$. Here the numerical value of $L(f_4,1)$ is
\begin{equation}
L(f_4,1)= 0.67496319716994177129269568273091339919322842904407 \cdots.
\end{equation}
The numerical value of the number $v^\perp$ is
\begin{equation}
v^\perp= 0.37369955695472976699767292752499463211766555651682 \cdots.
\end{equation}
The $j$-value of $\tau^\perp:=\frac{1}{2}+i \,v^\perp$ is a rational number
\begin{equation}
j(\tau^\perp)=\left(\frac{215}{28} \right)^3.
\end{equation}
LMFDB includes only one rationally defined elliptic curve with the above $j$-invariant, which has \textbf{14.2.a.a} as its eigenform. In fact, this elliptic curve is the modular curve $X_0(14)$
\begin{equation}
y^2+xy+y=x^3+4x-6.
\end{equation}
The readers are referred to \cite{Candelas} for more details.

\subsection{The computations of Deligne's periods}

Now, we are ready to compute the periods $c^\pm(h^3(\mathscr{X}_{-1/7})(2))$ for the critical motive $h^3(\mathscr{X}_{-1/7})(2)$. The matrix of the involution $F_\infty$ is given by the formula \ref{eq:complexconjugationinvolution}. Since the Picard-Fuchs operator \ref{eq:aesz34PF} does not have negative singularities, from Section \ref{sec:negativecase}, $F_\infty$ is constant in the interval $(-\infty, 0)$. Then formula \ref{eq:involutionnegativephi} and formula \ref{eq:prepotentialcoeff} tells us that $F_\infty$ is given by
\begin{equation}
F_\infty=
\left(
\begin{array}{cccc}
 1 & 1 & -6 & 12 \\
 0 & -1 & 12 & -24 \\
 0 & 0 & -1 & 0 \\
 0 & 0 & -1 & 1 \\
\end{array}
\right).
\end{equation}
The two linearly independent eigenvectors of $F_\infty$ associated to the eigenvalue $1$ are
\begin{equation}
v^+_1=(1,0,0,0),~v^+_2=(0,-12,0,1),
\end{equation}
hence the subspace of $H^3(X,\mathbb{Q})$ on which $F_\infty$ acts as $1$ is spanned by
\begin{equation}
\beta^0~\text{and}~-12\, \beta^1+\alpha^1.
\end{equation}
From Section \ref{sec:computationdeligneperiods}, the Deligne's period $c^+(h^3(\mathscr{X}_{-1/7})(2))$ is given by
\begin{equation}
c^+(h^3(\mathscr{X}_{-1/7})(2))=(2 \pi i)^4  \frac{1}{(2 \pi i)^6}  \det 
\left(
\begin{array}{cc}
\int_{X} \Omega_{-1/7} \smile \beta^0 &  \int_{X} \Omega_{-1/7} \smile (-12\,\beta^1+\alpha^1) \\
\int_{X} \Omega'_{-1/7} \smile \beta^0 & \int_{X} \Omega'_{-1/7} \smile  (-12\,\beta^1+\alpha^1) \\
\end{array}
\right),
\end{equation}
which is equal to
\begin{equation}
c^+(h^3(\mathscr{X}_{-1/7})(2))=12 \,\lambda^2 (2 \pi i)^4 \det 
\left(
\begin{array}{cc}
\varpi_0(-1/7) & -\varpi_1(-1/7)+\varpi_2(-1/7)\\
\varpi'_0(-1/7) & -\varpi'_1(-1/7)+\varpi'_2(-1/7) \\
\end{array}
\right).
\end{equation}
Since Deligne's period is only well-defined up to a nonzero rational multiple, we have the freedom to let $c^+(h^3(\mathscr{X}_{-1/7})(2))$ be
\begin{equation}
c^+(h^3(\mathscr{X}_{-1/7})(2))= \pi^4 \det 
\left(
\begin{array}{cc}
\varpi_0(-1/7), & -\varpi_1(-1/7)+\varpi_2(-1/7)\\
\varpi'_0(-1/7), & -\varpi'_1(-1/7)+\varpi'_2(-1/7) \\
\end{array}
\right).
\end{equation}
Plug in the numerical values of $\varpi^{(n)}_i(-1/7)$, we find that
\begin{equation} \label{eq:delignecaesz}
c^+(h^3(\mathscr{X}_{-1/7})(2))=-\frac{2401}{32} L(f_2,1)L(f_4,2)=-\frac{2401}{32} L(h^3(\mathscr{X}_{-1/7})(2),0),
\end{equation} 
which indeed satisfies the prediction of Deligne's conjecture.

Now, let us look at the period $c^-(h^3(\mathscr{X}_{-1/7})$. The two linearly independent eigenvectors of $F_\infty$ associated to the eigenvalue $-1$ are
\begin{equation}
v^-_1=(0,0,2,1),~v^-_2=(-1,2,0,0),
\end{equation}
hence the subspace of $H^3(X,\mathbb{Q})$ on which $F_\infty$ acts as $-1$ is spanned by
\begin{equation}
2 \alpha^0+\alpha^1~\text{and}~-\beta^0+2 \beta^1.
\end{equation}
From Section \ref{sec:negativecase}, the period $c^-(h^3(\mathscr{X}_{-1/7})$ is given by
\begin{equation}
c^-(h^3(\mathscr{X}_{-1/7})=
 \det 
\left(
\begin{array}{cc}
\left(32 \frac{\zeta(3)}{(2 \pi i)^3}-1\right)\varpi_0-2\varpi_1+12\varpi_2-8\varpi_3,& -\varpi_0+2\varpi_1\\
\left(32 \frac{\zeta(3)}{(2 \pi i)^3}-1\right)\varpi'_0-2\varpi'_1+12\varpi'_2-8\varpi'_3,& -\varpi'_0+2\varpi'_1 \\
\end{array}
\right),
\end{equation}
Notice that here we have thrown away a nonzero rational constant. Plug in the numerical values of $\varpi^{(n)}_i(-1/7)$, we find that
\begin{equation} \label{eq:deligneminusperiods}
c^-(h^3(\mathscr{X}_{-1/7}))=\frac{1029}{32}\,\pi^{-3} \frac{L(f_4,1)L(f_2,1)}{v^\perp}.
\end{equation}
A detailed study of this equation is presented in the paper \cite{Yang1}.

\section{Conclusion and further prospects} \label{sec:conclusion}

In this paper, we first briefly review the concept of pure motives, which plays a very crucial role in modern number theory and algebraic geometry. We try to illustrate the idea of pure motives through their classical realizations, which is perhaps easier to understand for physicists. Then we briefly discuss the construction of the $L$-function associated to a pure motive, and the (conjectured) analytic properties of $L$-functions. Next, we introduce Deligne's conjecture on the special values of the $L$-functions of critical motives.

As Deligne's conjecture is potentially extremely difficult to prove, therefore it is very interesting to see whether researches in other areas, e.g. string theory and mirror symmetry, could provide any insights into the conjecture itself. This is exactly the motivation of this paper! We have shown that mirror symmetry provides all the geometric data needed in the computation of the Deligne's period of a Calabi-Yau threefold. More precisely, we have developed a method to compute the Deligne's period of a smooth Calabi-Yau threefold in a one-parameter mirror family. We also illustrate how this method works by computing the Deligne's periods for two important examples.

In order to verify whether a Calabi-Yau threefold satisfies Deligne's conjecture, we also need to find the special value of its $L$-function at a critical integral point, which in practice is very difficult. In the paper \cite{Candelas}, the authors are able to find the $L$-function of a special Calabi-Yau threefold, and they have numerically computed the special values of this $L$-function. In this paper, we have numerically evaluated the Deligne's period for this special Calabi-Yau threefold, and we have numerically shown that it indeed satisfies Deligne's conjecture.

The results of this paper raise many interesting questions. For example, the computations in Sections \ref{sec:computationdeligneperiods} and \ref{sec:explicitDeligneperiods} are on the complex side, it is very interesting to ask whether they have interpretations on the K\"ahler side. More concretely, under the mirror map, do the Deligne's periods $c^\pm(h^3(\mathscr{X}_\varphi))$ have any interesting interpretations on the K\"ahler side? Another equally interesting question is whether the Deligne's periods $c^\pm(h^3(\mathscr{X}_\varphi))$ have any interesting interpretations in string theory or other related physics theories. Answering these questions might shed new lights on the nature of the Deligne's periods of Calabi-Yau threefolds, or even on the proof of Deligne's conjecture.

\section*{Acknowledgments}

The author is grateful to the referee of Nuclear Physics, Section B for pointing out several typos and raising several interesting questions.

\end{document}